\documentclass[aps,nofootinbibt, onecolumn, peprint]{revtex4}
\usepackage{eurosym}
\usepackage{graphicx}
\usepackage{amsmath}
\usepackage{bm}
\usepackage{color}
\usepackage{amsfonts}

\def\be{\begin{equation}}
\def\ee{\end{equation}}
\def\bea{\begin{eqnarray}}
\def\eea{\end{eqnarray}}

\begin{document}

\preprint{gr-qc/yymmnnn}
\title{Approximate Analytical Solution for the Dynamic Model of Large
Amplitude Non-Linear Oscillations Arising in Structural Engineering }
\author{J. A. Belinch\'{o}n}
\email{abelcal@ciccp.es}
\affiliation{Departamento de F\'{\i}sica Atomica, Molecular y Nuclear, Universidad
Computense de Madrid, E-28040 Madrid, Espa\~{n}a, and Universidad de Atacama, Facultad de Ciencias Naturales, Departamento de F\'isica, Copiap\'o, Chile}
\author{T. Harko}
\email{t.harko@ucl.ac.uk}
\affiliation{ Department of Physics, Babes-Bolyai University, Kogalniceanu Street,
Cluj-Napoca 400084, Romania,}
\affiliation{Department of Mathematics, University College London, Gower Street, London
WC1E 6BT, United Kingdom}
\author{M. K. Mak}
\email{mankwongmak@gmail.com}
\affiliation{Universidad de Atacama, Facultad de Ciencias Naturales, Departamento de F\'isica, Copiap\'o, Chile}
\date{\today }

\begin{abstract}
In this work we obtain an approximate solution of the strongly nonlinear
second order differential equation $\frac{d^{2}u}{dt^{2}}+\omega
^{2}u+\alpha u^{2}\frac{d^{2}u}{dt^{2}}+\alpha u\left( \frac{du}{dt}\right)
^{2}+\beta \omega ^{2}u^{3}=0$, describing the large amplitude free
vibrations of a uniform cantilever beam, by using a method based on the
Laplace transform, and the convolution theorem. By reformulating the initial
differential equation as an integral equation, with the use of an iterative
procedure, an approximate solution of the nonlinear vibration equation can
be obtained in any order of approximation. The iterative approximate
solutions are compared with the exact numerical solution of the vibration
equation.


\textbf{MSC classification codes: 74Hxx, 	41Axx, 65D99}

\end{abstract}

\maketitle


\section{Introduction}

Structural engineering is an important field dealing with the analysis and
design of structures that support or resist loads. It is commonly known that
the history of structural engineering contains many collapses and failures.
In order to avoid this kind of tragedy, structural engineers, and
researchers have attempted to develop various methodologies for handling
different kind of structures, for instance, building and bridge structures
according to physical laws and empirical knowledge of the structural
performance of various geometries and materials.

It is worth to note that many engineering structures can be modelled as a
slender, flexible cantilever beam, carrying a lumped mass with rotary
inertia at an intermediate point along its span \cite{Ham,a4}. The simplest
theoretical approach to these kind of systems is represented by the linear
analysis of vibration, where one introduces the approximation that the
frequency of free vibration of beam systems is independent of the amplitude
of the oscillation. The linear approximation is valid only when the
amplitude is relatively small. Since realistic beam systems are slender and
flexible, in solving practical engineering problems one should take into
account that they undergo relatively large amplitude flexural vibrations.
Due to their mathematical complexity, resulting from the highly nonlinear
structure, and the complicate character of the large amplitude oscillations,
usually this kind of nonlinear mathematical problems do not allow an exact
treatment, and therefore approximate techniques and perturbation methods
must be used for their study. In \cite{Ham} the single-term harmonic balance
method and two terms harmonic balance method was used, to obtain the
approximate solution to the period of oscillation. The model introduced in
\cite{Ham} has attracted considerable attention, and several approximate or
semi-analytical mathematical methods have been introduced for the study of
the free large amplitude nonlinear oscillations of a slender, inextensible
cantilever beam. Such methods are the homotopy perturbation method \cite{p1}%
, modified Lindstedt-Poincare methods \cite{p2,p3}, and the optimal homotopy
asymptotic method \cite{p4}. In \cite{p5} six different analytical methods
are applied to solve the dynamic model of the large amplitude non-linear
oscillation equation. The study of the systems experiencing large amplitude
vibrations is an important field of study in both mathematics and
engineering \cite{a1,a3,a6,a7,a8}. Furthermore, different mathematical
techniques and/or numerical procedures for handling the nonlinear
oscillation problems and, more generally, of the strongly nonlinear
differential equations, have been tremendously developed in \cite%
{10,11,12,13,14,15,16,17,18}. In many situations the nonlinear oscillation
equation can be reduced to a Li\'{e}nard type equation, which allows
obtaining exact analytical representations of the solutions of nonlinear
oscillators \cite{21}. For a study of anharmonic oscillations see \cite{22}.

The Laplace transform method is a very powerful approach for solving linear
differential equations, as well as systems of linear differential equations.
However, it can be also used to obtain approximate iterative solutions of
nonlinear differential equations. By using the Laplace transform and the
convolution theorem the nonlinear ordinary differential equation can be
transformed into an integral equation, whose solution can be obtained by the
method of successive approximations. Note that the Laplace transformation
method has been intensively applied to the study of the time evolution of
relativistic dissipative universe and the vacuum solutions of the
gravitational field equations in the brane world model \cite{01,02,03,04}.

The purpose of the present paper is to study the nonlinear differential
equation describing the free vibration of a uniform cantilever beam by using
its equivalent formulation as an integral equation obtained via the Laplace
transform. The iterative solutions of the vibration equation are explicitly
obtained in the first two orders of approximation, and they are compared
with the exact numerical solution.

The present paper is organized as follows. In Section~\ref{sect2} we briefly
present the derivation of the equation describing the free vibrations of the
cantilever beams. The exact solution of the vibrations equation is obtained
in Section~\ref{sect3}. The Laplace transform method for nonlinear ordinary
differential equations is introduced in Section~\ref{sect4}. The iterative
procedure and the first two approximations of the solution of the nonlinear
vibrations equation are presented in Section~\ref{sect5}. We discuss and
conclude or results in Section~\ref{sect6}.

\section{Large amplitude free vibrations of a uniform cantilever beam}

\label{sect2}

In the following we adopt the non-linear oscillation model introduced in
\cite{Ham}, which we briefly describe in the following. The basic physical
model consists of a clamped beam at the base, and free at the tip. The beam
caries a lumped mass $M$ and rotary inertia $J$ located at an arbitrary
intermediate point $s=d$ along its span. Following (\cite{Ham}) we introduce
the simplifying assumption of a uniform beam, and denote the constant length
by $l$, with $m$ representing the mass per unit length. The thickness of
this conservative simple beam is assumed to be very small as compared to its
length, and therefore we can safely ignore the effects of the rotary inertia
and of the shearing deformation. Moreover, as argued in \cite{Ham}, the beam
can be taken as inextensible, and this assumption implies that the length of
the beam neutral axis remains constant during the motion. Then, the elastic
potential $V$, due to the bending of the beam, can be written as \cite{Ham}
\begin{equation}
V=\frac{EIl}{2}\int_{0}^{1}{R^{2}\left( \xi ,t\right) d\xi },
\end{equation}%
where $EI$ is the modulus of flexural rigidity, $\xi =s/l$ is the
dimensionless arc length, and $R(\xi ,t)$ is the radius of curvature of the
beam neutral axis. The kinetic energy $T$ of the beam can be represented as
\begin{equation}
T=\frac{ml}{2}\int_{0}^{1}{\left( \dot{x}^{2}+\dot{y}^{2}\right) d\xi }+%
\frac{1}{2}M\left. \left( \dot{x}^{2}+\dot{y}^{2}\right) \right\vert _{\eta
l}+\frac{1}{2}J\left. \dot{\theta}^{2}\right\vert _{\eta l},
\end{equation}%
where the overdot denotes the derivative of time, $\eta =d/l$ is the
dimensionless relative position parameter of the attached inertia element,
and $\theta $ is the slope of the elastica, which can be obtained as $\sin
\theta =dy/ds$. In order to obtain the discrete single mode, single
coordinate beam Lagrangian one expands the kinetic and potential energies
into power series, retaining nonlinear terms up to fourth order, and looking
for an approximate single mode solution of the form $y(\xi ,t)=\Phi (\xi
)u(t)$.

Hence one obtains in this approximation the beam Lagrangian as \cite{Ham}
\begin{equation}
L=\frac{ml}{2}\left( \alpha _{1}\dot{u}^{2}+\alpha _{3}\lambda ^{2}u^{2}\dot{%
u}^{2}-\alpha _{2}\beta ^{2}u^{2}-\alpha _{4}\beta ^{2}\lambda
^{2}u^{4}\right) ,  \label{L}
\end{equation}%
where $\lambda $, $\beta $ and $\alpha _{i}$, for $i=1,2,3,4$ are constants.
After a rescaling of the time coordinate and of the constants we obtain the
equation of motion corresponding to the Lagrangian (\ref{L}) in the form
\cite{Ham}, \cite{p1,p2,p3,p4,p5}
\begin{equation}
\ddot{u}+\frac{\alpha u}{1+\alpha u^{2}}\dot{u}^{2}+\omega ^{2}\left( \frac{u%
}{1+\alpha u^{2}}+\frac{\beta u^{3}}{1+\alpha u^{2}}\right) =0,  \label{eqe}
\end{equation}%
where $\omega $ is a constant. Eq.~(\ref{eqe}), describing the dynamic model
of large amplitude nonlinear oscillations arising in the structural
engineering can be written in the equivalent form of the following second
order ordinary differential equation
\begin{equation}
\frac{d^{2}u}{dt^{2}}+\omega ^{2}u+\alpha u^{2}\frac{d^{2}u}{dt^{2}}+\alpha
u\left( \frac{du}{dt}\right) ^{2}+\beta \omega ^{2}u^{3}=0,  \label{1}
\end{equation}%
with the initial conditions given by
\begin{equation}
u\left( 0\right) =A,\frac{du}{dt}\left( 0\right) =0,  \label{2}
\end{equation}%
where $A$ is the arbitrary constant.

From a physical point of view, the third and fourth terms in Eq.~(\ref{1})
represent inertia-type cubic non-linearity, arising from the assumption of
the beam inextensibility. The last term in Eq.~(\ref{1}) is a static-type
cubic non-linearity, associated with the potential energy stored in bending.
There are two model constants, $\alpha $ and $\beta $, resulting from the
discretization procedure. The numerical values of these constants must be
determined from physical considerations.

\section{The exact solution of the nonlinear beam oscillation equation}

\label{sect3}

Now in order to solve Eq. (\ref{1}) exactly, we introduce the arbitrary
function $y\left( u\right) $, defined as
\begin{equation}
y\left( u\right) =\left( \frac{du}{dt}\right) ^{2},\frac{dy(u)}{du}=2\frac{%
d^{2}u}{dt^{2}}.  \label{u}
\end{equation}%
Thus, we rewrite Eq.~(\ref{1}) as
\begin{equation}
\frac{dy}{du}+\frac{2\alpha u}{1+\alpha u^{2}}y=-\frac{2\omega ^{2}u\left(
1+\beta u^{2}\right) }{1+\alpha u^{2}},  \label{3}
\end{equation}%
with the general solution given by

\begin{equation}
y\left( u\right) =\frac{C_{0}-\omega ^{2}u^{2}\left( 1+\frac{\beta }{2}%
u^{2}\right) }{1+\alpha u^{2}},  \label{4}
\end{equation}%
where $C_{0}$ is the arbitrary constant of integration. With the help of the
relation $y\left( u\right) =\left( du/dt\right) ^{2}$, Eq.~(\ref{4}) can be
integrated to give
\begin{equation}
t_{\pm }\left( x\right) -t_{0\pm }=\pm \int \sqrt{\frac{1+\alpha u^{2}}{%
C_{0}-\omega ^{2}u^{2}\left( 1+\frac{\beta }{2}u^{2}\right) }}du,  \label{5}
\end{equation}%
where $t_{0\pm }$ are the arbitrary constants of integration. With the help
the initial conditions $u\left( 0\right) =A,\frac{du}{dt}\left( 0\right) =0$%
, we get the relation
\begin{equation}
C_{0}=\omega ^{2}A^{2}\left( 1+\frac{\beta }{2}A^{2}\right) .
\end{equation}

Thus we have obtained the exact solution for this specific nonlinear
oscillation problem, important in engineering applications. However, from a
practical point of view the integral representation of the solution as given
in Eq.~(\ref{5}) is not particularly useful. Therefore in the following we
will look for some approximate solutions of the nonlinear high amplitude
beam vibrations.

Obviously, the differential Eq.~ (\ref{1}) contains two linear terms $\frac{%
d^{2}u}{dt^{2}}$ and $\omega ^{2}u$, describing a simple harmonic
oscillation. By neglecting the other terms in Eq.~(\ref{1}), and by
integrating the corresponding linear differential equation, yields the zero
order approximate solution of Eq.~(\ref{1}). Starting from the linear
approximation we shall be able to present the approximate analytical
solution of Eq.~(\ref{1}) by using the Laplace transformation method in the
next Section.

\section{The Laplace transform method for nonlinear differential equations}

\label{sect4}

Even that the Laplace transform method is especially useful in solving the
Cauchy problem for linear differential equations, or systems of linear
differential equations, the method can sometimes be applied to find the
solution of the Cauchy problem for nonlinear differential equations. The
formal scheme for the application of the Laplace transform method is as
follows \cite{Teo}. Let's consider a non-linear differential equation of the
form
\begin{equation}
a_{n}y^{(n)}+a_{n-1}y^{(n-1)}+...+a_{1}y^{\prime }+a_{0}y+g\left[
t,y,y^{\prime (m)}\right] =f(t),  \label{Ln}
\end{equation}%
where $m\leq n$, $a_{j}\in \mathbb{C}$, $f:\Re _{+}\rightarrow \mathbb{C}$, $%
g:\Re _{+}\times \Re ^{m+1}\rightarrow \mathbb{C}$. The initial conditions
for Eq.~(\ref{Ln}) are
\begin{equation}
y(0)=y_{0},y^{\prime }(0)=y_{1},...,y^{(n-1)}(0)=y_{n-1}.
\end{equation}%
By applying the Laplace transform operator $L$ to Eq.~(\ref{Ln}) we obtain
an equality of the form
\begin{equation}
P_{1}(p)\cdot L(y)(p)-P_{2}\left( p,y_{0},y_{1},...,y_{n-1}\right) +L\left[
g\left( t,y,y^{\prime (m)}\right) \right] (p)=L(f)(p),  \label{L1}
\end{equation}%
where
\begin{equation}
P_{1}(p)=a_{n}p^{n}+a_{n-1}p^{n-1}+...+a_{1}p+a_{0},
\end{equation}%
and $P_{2}$ is a polynomial of the maximum degree $n-1$. Hence from Eq.~(\ref%
{L1}) we obtain the relation
\begin{equation}
L(y)(p)=G_{1}(p)H(p)-G_{1}(p)\cdot L\left[ g\left( t,y,y^{\prime (m)}\right) %
\right] (p),
\end{equation}%
where we have denoted
\begin{equation}
G_{1}=\frac{1}{P_{1}},H=L(f)+P_{2}.
\end{equation}

With the use of the inverse Laplace transform we obtain
\begin{eqnarray}  \label{L2}
y(t)&=&g_1(t)*h(t)-g_1(t)*g\left[t,y,y^{\prime (m)}\right]%
=\left(g_1*h\right)(t)-  \notag \\
&&\int_0^t{g_1(t-x)g\left[x,y(x),y^{\prime (m)}(x)\right]dx},
\end{eqnarray}
where
\begin{equation}
g_1=L^{-1}\left(G_1\right),h=L^{-1}(H).
\end{equation}

Eq.~(\ref{L2}) is a nonlinear integral equation for the unknown function $y$%
. If Eq.~(\ref{Ln}) does admit a solution, the solution can be obtained by
using the method of successive approximations. We take
\begin{equation}
y_{0}(t)=\left( g_{1}\ast h\right) (t),
\end{equation}%
and then the successive steps of the iterative solution can be constructed
with the help of the equations
\begin{equation}
y_{k+1}(t)=\left( g_{1}\ast h\right) (t)-\int_{0}^{t}{g_{1}(t-x)g\left[
x,y_{k}(x),y_{k}^{\prime }(x),...,y_{k}^{(m)}(x)\right] dx},k\in \mathbb{N}.
\end{equation}

The iterative procedure can be also applied to Eq.~(\ref{L1}) itself, which
may be more advantageous in some cases. Thus, by taking
\begin{equation}
Y_{0}(p)=\left( G_{1}\cdot H\right) (p),
\end{equation}%
from Eq.~(\ref{L1}) we obtain
\begin{equation}
Y_{k+1}(p)=\left( G_{1}\cdot H\right) (p)-G_{1}(p)\cdot L\left\{ g\left[
t,y_{k}(t),y_{k}^{\prime }(t),...,y_{k}^{(m)}(t)\right] \right\} (p),
\end{equation}%
where $y_{k}=L^{-1}\left( Y_{k}\right) $, $k\in \mathbb{N}$. If $Y_{0}$ is
known, we obtain immediately $y_{0}=L^{-1}\left( Y_{0}\right) $, which
allows us to determine
\begin{equation}
Y_{1}(p)=\left( G_{1}\cdot H\right) (p)-G_{1}(p)L\left[ g(t),y_{0},y_{0}^{%
\prime },...,y_{0}^{(m)}\right] (p),
\end{equation}%
leading to $y_{1}=L^{-1}\left( Y_{1}\right) $ etc.

\section{The Laplace transform method for the equation of the free
vibrations of the cantilever beam}

\label{sect5}

Due to the highly non-linear structure of the differential Eq. (\ref{1}), it
is difficult to find the exact solution of it, given by $u=u\left( t\right) $%
, in a form practical for explicit computations. In order to solve Eq.~(\ref%
{1}) directly, we apply the Laplace transformation and the convolution
theorem. After taking the Laplace transform of Eq.~(\ref{1}), we obtain
\begin{equation}
U(p)=\frac{Ap}{p^{2}+\omega ^{2}}-\frac{1}{p^{2}+\omega ^{2}}L[\alpha u^{2}%
\ddot{u}+\alpha u\dot{u}^{2}+\beta \omega ^{2}u^{3}](p),
\end{equation}%
where $U(p)=(L\cdot u)(p)$. With the use of the inverse Laplace transform
and of the convolution theorem, we obtain the formal representation of the
solution of Eq.~(\ref{1}) as
\begin{equation}
u(t)=A\cos (\omega t)-\frac{1}{\omega }\int_{0}^{t}\sin {\left[ \omega (t-x)%
\right] u(x)}\left[ \alpha u(x)u^{\prime \prime }(x)+\alpha u^{\prime
2}(x)+\beta \omega ^{2}u^{2}(x)\right] {dx}.
\end{equation}%
The approximate solutions of the above integral equation can be obtained
iteratively as
\begin{eqnarray}
u_{n}(t) &=&A\cos (\omega t)-  \notag \\
&&\frac{1}{\omega }\int_{0}^{t}\sin {\left[ \omega (t-x)\right] u_{n-1}(x)}%
\left[ \alpha u_{n-1}(x)u_{n-1}^{\prime \prime }(x)+\alpha u_{n-1}^{\prime
2}(x)+\beta \omega ^{2}u_{n-1}^{2}(x)\right] {dx},n\in \mathbb{N}.
\label{qq}
\end{eqnarray}%
Thus $u_{n}(t)$ gives the $n^{th}$ order approximate solution of Eq. (\ref{1}%
). Hence, through the method of successive approximations, we have obtained
the complete approximate analytical solution of Eq.~(\ref{1}) describing the
dynamic model of large amplitude non-linear oscillations arising in the
structural engineering, given by
\begin{equation}
u(t)=\lim_{n\rightarrow \infty }u_{n}(t).  \label{a3}
\end{equation}

In the first order of approximation, by taking $u_{0}(t)=A\cos (\omega t)$
we obtain
\begin{equation}
u_{1}(t)\approx A\cos (\omega t)+\frac{1}{16}A^{3}\sin (\omega t)\left[
2\omega t(2\alpha -3\beta )+(2\alpha -\beta )\sin (2\omega t)\right] .
\label{first}
\end{equation}%
In the second order of approximation we find
\begin{eqnarray}
u_{2}(t) &\approx &A\cos (\omega t)+  \notag  \label{sec} \\
&&\frac{A^{3}}{188743680}\Bigg\{18(2\alpha -\beta )^{3}(18\alpha -\beta
)\cos (9t\omega )A^{6}-120t(2\alpha -3\beta )(34\alpha -3\beta )(\beta
-2\alpha )^{2}\omega \times   \notag \\
&&\sin (7t\omega )A^{6}-5(\beta -2\alpha )^{2}\left[ 500A^{2}\alpha
^{2}+\left( 6528-420A^{2}\beta \right) \alpha +9\beta \left( 9A^{2}\beta
-64\right) \right] \cos (7t\omega )\times   \notag \\
&&A^{4}+240t\left[ 4\alpha ^{2}-8\beta \alpha +3\beta ^{2}\right] \left[
39A^{2}\beta ^{2}-12\left( 11\alpha A^{2}+16\right) \beta +4\alpha \left(
35\alpha A^{2}+288\right) \right] \omega \times   \notag \\
&&\sin (5t\omega )A^{4}-20(2\alpha -\beta )\Bigg[-488\alpha
^{3}A^{4}+477\beta ^{3}A^{4}-1422\alpha \beta ^{2}A^{4}+144t^{2}(2\alpha
-3\beta )^{2}\times   \notag \\
&&(6\alpha -\beta )\omega ^{2}A^{4}+1500\alpha ^{2}\beta A^{4}-4992\alpha
^{2}A^{2}-3168\beta ^{2}A^{2}+7296\alpha \beta A^{2}-55296\alpha +9216\beta %
\Bigg]\times   \notag \\
&&\cos (5t\omega )A^{2}+720t(2\alpha -3\beta )\omega \Bigg[-344\alpha
^{3}A^{4}+543\beta ^{3}A^{4}-1218\alpha \beta ^{2}A^{4}+16t^{2}(2\alpha
-3\beta )^{2}\times   \notag \\
&&(2\alpha -\beta )\omega ^{2}A^{4}+1076\alpha ^{2}\beta A^{4}-1024\alpha
^{2}A^{2}-1920\beta ^{2}A^{2}+2304\alpha \beta A^{2}-6144\alpha +3072\beta %
\Bigg]\times   \notag \\
&&\sin (3t\omega )A^{2}+\Bigg[215408\alpha ^{4}A^{6}+635787\beta
^{4}A^{6}-11520t^{4}(2\alpha -3\beta )^{4}\omega ^{4}A^{6}-1826328\alpha
\beta ^{3}A^{6}+  \notag \\
&&2131000\alpha ^{2}\beta ^{2}A^{6}-1136928\alpha ^{3}\beta
A^{6}-299520\alpha ^{3}A^{4}-2376000\beta ^{3}A^{4}+3653760\alpha \beta
^{2}A^{4}-  \notag \\
&&1340160\alpha ^{2}\beta A^{4}-6635520\alpha ^{2}A^{2}+4239360\beta
^{2}A^{2}+1440t^{2}(2\alpha -3\beta )\times   \notag \\
&&\Bigg(-705\beta ^{3}A^{4}+2\left( 673\alpha A^{2}+384\right) \beta
^{2}A^{2}+8\alpha \left( A^{2}\alpha \left( 15A^{2}\alpha -64\right)
-256\right) -  \notag \\
&&4\left( 9A^{2}\alpha -16\right) \left( 23\alpha A^{2}+48\right) \beta %
\Bigg)\omega ^{2}A^{2}+737280\alpha \beta A^{2}+11796480\alpha -5898240\beta %
\Bigg]\cos (t\omega )+  \notag \\
&&180\Bigg[-1264\alpha ^{4}A^{6}-3583\beta ^{4}A^{6}+10392\alpha \beta
^{3}A^{6}-12264\alpha ^{2}\beta ^{2}A^{6}+6624\alpha ^{3}\beta
A^{6}+1280\alpha ^{3}A^{4}+  \notag \\
&&13536\beta ^{3}A^{4}-21568\alpha \beta ^{2}A^{4}+64t^{2}(2\alpha -3\beta
)^{2}\Bigg(A^{2}\alpha ^{2}+\left( 48-17A^{2}\beta \right) \alpha +3\beta
\left( 3A^{2}\beta -8\right) \Bigg)  \notag \\
&&\omega ^{2}A^{4}+8832\alpha ^{2}\beta A^{4}+24576\alpha
^{2}A^{2}-24576\beta ^{2}A^{2}+4096\alpha \beta A^{2}-65536\alpha
+32768\beta \Bigg]\times   \notag \\
&&\cos (3t\omega )+960t\omega \Bigg[16t^{2}(2\alpha -3\beta )^{2}\left(
2A^{2}\alpha ^{2}+\left( 16-11A^{2}\beta \right) \alpha +3\beta \left(
5A^{2}\beta -8\right) \right) \omega ^{2}A^{4}+  \notag \\
&&3\Bigg(-1056\beta ^{4}A^{6}+\left( 2807\alpha A^{2}+1080\right) \beta
^{3}A^{4}-2\left( A^{2}\alpha \left( 1449\alpha A^{2}+1336\right)
-3072\right) \beta ^{2}A^{2}-  \notag \\
&&8\alpha \left( 3A^{2}\alpha \left( 3\alpha \left( 3\alpha A^{2}+8\right)
A^{2}+256\right) -2048\right) +4\Bigg(A^{2}\alpha \left( \alpha \left(
333\alpha A^{2}+584\right) A^{2}+768\right) -  \notag \\
&&6144\Bigg)\beta \Bigg)\Bigg]\sin (t\omega )\Bigg\}.
\end{eqnarray}

\begin{figure}[h]
\centering
\includegraphics[width=8cm]{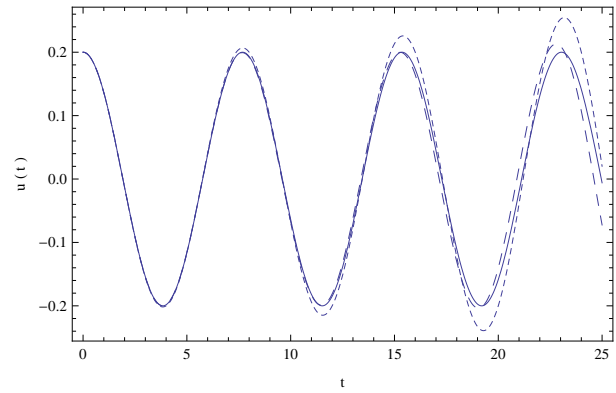}
\caption{Time variation of $u$ for $\protect\alpha =1.2$, $\protect\beta %
=3.7 $, $\protect\omega = \protect\pi/4$ and $A=0.20$. The solid curve
represents the exact numerical solution, the dotted curve depicts the first
order Laplace transform approximation, while the dashed curve gives the
solution in the second Laplace transform approximation.}
\label{fig1}
\end{figure}

In Fig.~\ref{fig1} we present the time variation of $u$, obtained by solving
numerically Eq.~(\ref{1}), in the first order approximation, with the use of
Eq.~(\ref{first}), and in the second order approximation, given by Eq.~(\ref%
{sec}).

\section{Conclusions}

\label{sect6}

In the present paper we have considered the possibility of using the Laplace
transform and the convolution theorem for solving some nonlinear
differential equations. With the use of the Laplace transform nonlinear
differential equations can be easily reformulated as integral equations,
thus opening the way towards obtaining iterative solutions. As an
application of this method we have considered the case of the equation
describing the free non-linear vibrations of a uniform cantilever beam. We
have reformulated this equation as an integral equation, by using the
Laplace transform method, and we have explicitly obtained the first two
terms in the iterative expansion of the solution. The comparison with three
exact numerical solutions shows that the obtained expressions give a good
approximation of the exact solution, at least for small numerical values of $%
A$. Already the first order approximation, which has a very simple form, can
provide a high precision approximation.

Therefore the Laplace transform method of solution of non-linear
differential equations may provide a powerful tool for the investigation of
the properties of the highly non-linear ordinary differential equations.


\end{document}